\documentclass{article}
\usepackage{amssymb, amsfonts, amsmath, amsthm}

\usepackage[russian]{babel}
\newtheorem{theorem}{Теорема}

\newtheorem{cor}{Следствие}

\newcommand{\ep}{\epsilon}

\newcommand{\R}{\mathbb R}   

\begin{document}
\title{Замечания о чебышевских координатах}

\author{ Ю.~Бураго, С.~Иванов, С.~Малев}
\thanks{Исследования первых двух авторов частично поддержаны грантом РФФИ 05-01-00939}

\date{}

\maketitle

Под поверхностью всюду, если не оговорено противное, мы  имеем в виду
двумерное риманово многообразие.
Как известно, система координат на поверхности $M$ называется
чебышевской, если первая основная форма поверхности имеет в этих координатах
вид

$$ds^2=du^2+2\cos \theta(u,v) dudv+dv^2.$$

Иными словами, все координатные линии параметризованы длиной дуги,
а $ \theta(u,v)$ --- угол между координатными линиями.

Формула Хацидакиса (см., например, \cite[стр. 413]{Bak_Ver_Kant}
показывает, что на полной односвязной поверхности,
для которой $|\iint KdS|>2\pi$,
единые  чебышевские координаты могут не существовать.
Здесь $K$ --- гауссова кривизна, а
$S$ --- площадь.

Нас будет интересовать условия, при которых такие
координаты заведомо существуют. Наряду с римановыми 2-многообразиями мы будем
рассматривать более общий класс {\em поверхностей Александрова} (``называемых также двумерные
многообразия ограниченной кривизны'').
Этот  класс включает  в себя как
2-многообразия с римановыми или  полиэдральными ме\-три\-ка\-ми, так и их равномерные
пределы при условии равномерной ограниченности кривизны и сохранения структуры
многообразия, подробнее см. \cite{AlexZalg}
или \cite{Reshetnyak}.

П.~Л.~Чебышевым было доказано существование локальных чебышевских координат
в традиционном в то время предположении аналитичности.
Первый нелокальный результат  принадлежит И.~Бакельману \cite{Bakelman},
который доказал, что на  поверхности Александрова можно построить
единую чебышевскую сеть во всем  секторе между двумя ортогональными геодезическими,
если только положительная и отрицательная части интегральной кривизны этого сектора
обе меньше $\frac{\pi}{2}$. Нетрудно видеть, что для сектора этот
результат неулучшаем.

Много позже сходный
результат, но только в  случае $C^\infty$-гладких римановых многообразий, был получен
С.~Самелсон и В.~Даяванза  \cite{SmlsDyw}:
при тех же условиях на кривизну они
доказали существование единой чебышевской сети  на полной односвязной
$C^\infty$-гладкой поверхности. Их доказательство в основной своей части ---
аналитическое. Как известно, вопрос о локальном существовании чебышевских координат
сводится к решению определенной гиперболической квазилинейной системы уравнений
второго порядка; в \cite{SmlsDyw} эта система решена методом работы \cite{LadShub}.
\medskip

Мы покажем, что из результата  И.~Бакельмана следует такая теорема.

\begin{theorem}
\label{thm1}
Пусть $M$ --- полное односвязная поверхность Александрова,
удовлетворяющая условиям $\omega^+(M)<2\pi, \;\omega^-(M)<2\pi$.
Тогда на $M$ су\-ще\-ству\-ет единая
чебышевская система координат.

Если же  $\omega^+(M)<2\pi-\ep, \;\omega^-(M)<2\pi-\ep$, где $\ep>0$,
то существует такая чебышевская сеть, что углы между координатными линиями
равномерно отделены от  0 и $\pi$ (например, числом
$\ep/4$).
\end{theorem}

Здесь $\omega^+$     и $\omega^-$    --- положительная и отрицательная
кривизны $M$, рассматриваемые как
меры; в случае римановой метрики

\centerline{$\omega^+(M)=\iint_M K^+dS,\quad \omega^-(M)=\iint_M K^-dS$,}

\noindent где
 $a^+=\max\{a, 0\},\;a^-=\max\{-a, 0\}$.

 Как было отмечено, условия теоремы являются неулучшаемыми.

Для доказательства было бы достаточно разбить $M$ двумя пересекающимися
геодезическими на четыре сектора так, чтобы
для каждого из них выполнялись условия теоремы И.~Бакельмана.
Однако такое разбиение не всегда возможно и конструкцию приходится
несколько усложнить.

На самом деле рассуждения   И.~Бакельмана  доказывают несколько
более общее утверждение, чем было упомянуто выше. Именно,
рассмотрим на полной односвязной не\-огра\-ни\-чен\-ной
поверхности Александрова сектор $Q$, ограниченной кривыми
$\gamma_1, \;\gamma_2$ с общим началом $p=\gamma_1(0)=\gamma_2(0)$
в вершине сектора. Обозначим через $\tau$ поворот с внутренней
стороны границы сектора $Q$, за вычетом его вершины. Поворот
рассматривается как функция множества (заряд); через
$\tau^+,\;\tau^-$ обозначим положительную и отрицательную части
заряда $\tau$. (В случае гладкой кривой $\gamma\colon (0,a)\to M$
на римановом 2-многообразии $M$ \; $\tau^\pm (E)=\int_E k_g^\pm
ds$, где  $k_g$ --- геодезическая кривизна со стороны $Q$, а
интегрирование (по длине дуги) ведется по множеству $E\subset
(0,a)$.) Наконец, положим
$\tilde\omega^\pm(M)=\omega^\pm(M)+\tau^\pm(\gamma_1)+\tau^\pm(\gamma_2).$

\begin{theorem}
\label{Bak}
{\em (И.~Бакельман)}
Если для сектора $Q$ выполнены условия:
$$
\tilde\omega^+(Q)<\alpha, \quad \tilde\omega^-(Q)<\pi-\alpha,
$$
где $\alpha$ --- угол сектора $Q$ в его вершине $p$,
то в  $Q$ существует такая глобальная
чебышевская система координат, при которой кривые $\gamma_1, \;\gamma_2$
являются координатными линиями.

При этом сетевые углы отделены от 0 и $\pi$    минимумом величин
$\alpha-\tilde\omega^+(Q)$ и $\pi-\alpha-\tilde\omega^-(Q)$.
\end{theorem}

Ясно, что для доказательства теоремы \ref{thm1}
достаточно разбить $M$
на четыре сектора, удовлетворяющие условиям теоремы \ref{Bak}.
Такое разбиение достаточно осуществить только в случае
римановой  или  полиэдральной метрики, так как после этого общий случай
получается стандартным предельным переходом.
\medskip

Проблема существования чебышевских координат тесно связана с вопросом об
условиях, при которых между двумя поверхностями существуют билипшицевы
отображения  с контролируемыми константами Липшица,
см. \cite{BelBur, Bur}. Грубо говоря, постоянная Липшица оценивается через сетевые углы,
т.~е. через $\min \{\inf\theta,\,\inf(\pi-\theta)\}$.

Хотя при $\omega^-(M)>2\pi$ поверхность $M$ может не допускать чебышевских координат,
на ней все же можно ввести ``чебышевские координаты'' с ветвлением в конечном
числе точек, и такие обобщенные координаты
обеспечивают билипшицевое отображение на плоскость с контролируемой,
постоянной Липшица.
Например, для поверхностей, гомеоморфных плоскости, М.~Бонком и У.~Лангом
в \cite{BonkLang} доказано, что
при условии $\omega^+(M)<2\pi-\ep<2\pi, \;\omega^-(M)<C<\infty$
поверхность $M$ билипшицево эквивалентна плоскости с константой Липшица
$ \ep^{-\frac12}(2\pi+C)^{\frac12}$. Из существования
``чебышевских координаты с ветвлениями'' немедленно следует ослабленная
(т.~е. с худшей постоянной) теорема Бонка--Ланга.

Доказательство существования чебышевских координат с ветвлениями несущественно отличается
построений И.~Бакельмана; не описывая его конструкцию в деталях,
упомянем только, что сеть строится сначала на поверхности с полиэдральной метрикой;
общий случай получается предельным переходом. В условиях теоремы \ref{Bak},
сектор с полиэдральной метрикой
можно разбить на плоские параллелограммы, прилегающие по
целым сторонам и образующие решетку из строк и столбцов. Для построения координат
каждый параллелограмм заменяется прямоугольным параллелограммом на  $\R^2$ со сторонами
той же длины, так что новые параллелограммы замощают сектор в евклидовой плоскости.
В случае слишком большой отрицательной кривизны разбить сектор
на параллелограммы, образующие решетку,
может оказаться невозможно. Однако можно изменить конструкцию, допустив
примыкание более, чем четырех параллелограммов к некоторым вершинам в секторе (таких
вершин будет конечное число) и отображая часть параллелограммов на уже не обязательно
прямоугольные параллелограммы евклидовой плоскости.
 \medskip

Пусть  $M$ --- полное односвязное риманово 2-многообразие.
Назовем {\em крестом} фигуру на $M$, образованную тремя
такими простыми
натурально параметризованными геодезическими $\gamma_i$, что

(i) $\gamma_1\colon (-\infty, \infty)\to M$,
 \; $\gamma_i\colon (0, \infty)\to M$ при $i=2,3$;

(ii) $\gamma_1$  разбивает $M$ на две области (полуплоскости);

(iii) $\gamma_2$ и
 $\gamma_3$ начинаются на  $\gamma_1$
 и проходят в разных полуплоскостях относительно  $\gamma_1$, разбивая каждую
 полуплоскость на два сектора.

Тем самым
эти геодезические разбивают $M$ на четыре сектора $Q_i$. Обозначим через
$\alpha_i$ угол сектора $Q_i$ в его вершине $O_i$.

Аналогичную фигуру, но образованную квазигеодезическими на поверхности Александрова,
будем называть {\em обобщенным крестом}.

\begin{theorem}
\label{cross}
Пусть M - полное двумерное риманово многообразие, гомеоморфное плоскости и
удовлетворяющее условиям
\begin{equation}
\label{two_pi}
\omega^+ < 2\pi - 4\ep, \quad \omega^-<2\pi-4\ep
\end{equation}
 для
некоторого $\ep>0$. Тогда на $M$ существует такой крест,
что при всех $i=1,2,3,4$ \;
$\omega^+(Q_i)\leq\alpha_i-\ep$,\;
$\omega^-(Q_i)\leq\pi-\alpha_i-\ep$ и при этом для любых
двух точек $\gamma_i(t),\;\gamma_i(t')$, $i=1,2,3$, выполняется
\begin{equation}
\label{chord}
|t-t'|\sin\ep\leq d_M(\gamma_i(t),\gamma_i(t')).
\end{equation}
\end{theorem}

Аппроксимируя поверхность Александрова римановыми многообразиями, мы немедленно
получаем такое следствие.

\begin{cor} Пусть полная гомеоморфная плоскости
поверхность Александрова $M$
удовлетворяет условию \eqref{two_pi}. Тогда на $M$  существует обобщенный крест,
обладающий теми же свойствами, что крест в теореме \ref{cross}, но с заменой
$\omega^\pm$   на $\tilde\omega^\pm$, соответственно.
\end{cor}

В свою очередь, из этого следствия и теоремы \ref{Bak} немедленно
вытекает теорема \ref{thm1}.

{\em Замечание.} Теорема \ref{thm1} не утверждает, что в случае достаточно
гладкой римановой метрики построенные чебышевские координаты также будут гладкими.
Использование
аналитического метода  Самелсон и Даяванза  \cite{SmlsDyw} вместо конструкции
Бакельмана позволяет
гарантировать гладкость координатной сети в каждом секторе, но на границе секторов
координатные линии могут иметь изломы. Вопрос о существовании
 всюду гладкой  чебышевской сети остается пока открытым.

\medskip
 Доказательство  теоремы \ref{cross} следует тому же плану,
 что доказательство предложения 6.1
 работы \cite{BonkLang}.
\medskip

{\em Доказательство теоремы \ref{cross}.}

{\em Шаг 1.} Согласно только что  упомянутому предложению из
работы  Бонка--Ланга, существует геодезическая $\gamma_1$,
удовлетворяющая неравенству \eqref{chord}  и разбивающая $M$ на
две полуплоскости $M_k$, $k=1,2$, несущие равные порции
положительной и отрицательной кривизн.
Таким образом,  для каждой полуплоскости имеем
 \begin{equation}
 \label{curv}
 \omega^+(M_k)<\pi - 2\ep, \quad
 \omega^-(M_k)<\pi - 2\ep.
 \end{equation}
Рассмотрим для определенности $M_1$.  Нам надо найти в  $M_1$
геодезическую с началом на $\gamma_1$, разбивающий $M_1$ на два
сектора в соответствии с утверждением теоремы. Сразу отметим, что
неравенство \eqref{chord}  следует непосредственно из первого
неравенства \eqref{curv} и известного неравенства между дугой и
хордой, см. \cite{AlexZalg}, параграф 2 главы IX.

Как и в  \cite{BonkLang}, можно считать, что наша метрика плоская
вне некоторого компакта $\tilde D\subset\text{int}\,M_1$,
причем этот компакт гомеоморфен замкнутому диску, ограничен гладкой кривой,
тотально выпуклый и не пересекается с
$\gamma_1$. Кроме того, можно предполагать, что полуплоскость $M_1$ также тотально
выпукла; для этого достаточно, например,  отрезать  $M_2$ от  $M_1$ и
подклеить к $M_1$ евклидову полуплоскость вместо $M_2$.

{\em Шаг 2.} Обозначим через $S\tilde D$ расслоение единичных окружностей над
$\tilde D$. (В наших условиях это --- просто прямое произведение $\tilde D$ на
окружность $S^1$.)
Зададим следующее отображение $\sigma$ из $S\tilde D$ в
множество (максимально продолженных) геодезических поверхности $M_1$.
 Если $v=(p,\theta)\in S\tilde D$, то $\sigma(v)$ -  геодезическая,
 проходящая через точку $p$ ортогонально к $v$, причем ориентированная
 так, что вектор скорости
 геодезической $\sigma(v)$
 в точке $p$  и вектор $v$  определяют положительную ориентацию   $M_1$.

В силу теоремы Гаусса--Бонне и наших условий, геодезическая $\sigma(v)$
не имеет самопересечений и
разбивает $M_1$ на две области, из которых одна может быть ограниченной.
Обозначим через  $U(v)$ ту из областей, в которую направлен вектор $v$.
  Теперь введем на $S\tilde D$ угловую функцию
$\alpha(v)$ следующим образом: если  $U(v)$ --- полоса,
то $\alpha(v)=0$;  если $U(v)$ ---
полуплоскость, ограниченная только геодезической $\sigma(v)$,  то
$\alpha(v)=\pi$;  если $U(v)$ --- сектор,
то $\alpha(v)$ --- угол этого сектора; если  $U(v)$ ---
двуугольная область и оба угла этой области $\alpha_1$ и $\alpha_2$ ---
острые, то $\displaystyle{\alpha(v)=\frac{\alpha_1^2+\alpha_2^2}{\alpha_1+\alpha_2}}$;
если же один из углов (например, $\alpha_1$) двуугольника не острый, то
полагаем
$\displaystyle{\alpha(v)=\frac{\pi^2/4+\alpha_2^2}{\pi/2+\alpha_2}+\alpha_1-\frac{\pi}{2}}$;
 наконец, если $U(v)$ --- дополнение  $M_1$ до двуугольника, то $\alpha(v)=\pi-\alpha(-v)$.
Таким образом, мы определили угловую функцию $\alpha$ на всем   $S\tilde D$.
Сразу заметим, что в силу теоремы Гаусса--Бонне не острым может быть не более, чем
один угол двуугольника.

Нетрудно проверить, что функция $\alpha$ непрерывна и обладает свойством
 $$\alpha(v)+\alpha(-v)=\pi .$$

Нормируем меры $\omega^+$ и $\omega^-$, введя две новые меры:

$$W^{\pm}(A)=(\pi-2\ep) \frac{\omega^{\pm}(A)}{\omega^{\pm}({\tilde D})}.$$

Рассмотрим отображение $\phi: S\tilde D \rightarrow R^2$, которое
каждому вектору $v\in S\tilde D$ ставит в
соответствие точку
$$(W^+(U(v))-\alpha+\ep,\quad W^-(U(v))-(\pi-\alpha)+\ep).$$
Заметим, что $\phi(-v)=-\phi(v)$, в частности, если
$\phi(v)=(0,0)$, то также  $\phi(-v)=(0,0)$.

Для доказательства теоремы достаточно
найти такой вектор $v_0$, что $\phi(v_0)=(0,0)$.

\medskip
Зафиксируем любую точку $O$  границы диска $\tilde D$. Выберем
 в пространстве $S\tilde D$ две петли, $\gamma_1$ и $\gamma_2$, с вершиной $O$.

Петля  $\gamma_1$ представляет собой вращение вектора $v$ вокруг
точки $O$, начиная с вектора, перпендикулярного
к границе $\tilde S$ диска $\tilde D$ и направленного внутрь этого диска.

Петля  $\gamma_2$ получена обходом  $\tilde S$ (начиная с точки $O$) при условии,
что в каждой точке берется вектор $v$, перпендикулярный к $\tilde S$
и направленный внутрь $\tilde D$.

Эти петли очевидно  гомотопны. Действительно, достаточно
рассмотреть такой  диффеоморфизм $g$  диска
 $\tilde D$   на  евклидов единичный круг $D$, что $dg$ конформно
на границе  $\tilde D$, и, считая    $O'=g(O)$ нуль-вектором, взять $dg$-прообразы
семейства окружностей $t\,\partial D$, $0\leq t\leq 1$,
 вместе с полями их  нормалей, направленных внутрь.

Отображение $\phi$ переводит петли  $\gamma_1$, $\gamma_2$ в
петли пространства $\text{Im}\,\phi$; значит, последние гомотопны между собой в
 $\text{Im}\,\phi$.

Так как образ петли  $\gamma_2$ лежит на прямой $x+y=\pi-2\ep$, то он стягиваем в
 $\text{Im}\,\phi$.

В то же время, по выбору угловой функции $\alpha$, петля  $\gamma_1$ обладает свойством
$\phi(v)+\phi(-v) = (0,0)$. Отсюда следует, что индекс точки $(0,0)$
относительно петли  $\gamma_1$ не равен нулю и тем самым
она не стягиваема в $\R^2\setminus (0,0)$, а тем более в  $\text{Im}\,\phi$.
Действительно, угол, на который
поворачивается вектор с началом $(0,0)$ при прохождении петли  $\gamma_1$, есть сумма двух
равных друг другу углов вида $(2k+1)\pi$.

Таким образом, если нуль-вектор плоскости  не лежит в $\text{Im}\,\phi$
(как мы предположили), то петля  $\gamma_1$
не стягиваема, а вторая стягиваема. Это противоречит тому, что  петли
гомотопны.

Итак, найдется вектор $v$ единичного расслоения, для которого $\phi(v)= (0,0)$.
Иными словами, мы имеем
$W^+(U(v))=\alpha-\ep, \quad W^-(U(v)=\pi-\alpha-\ep$.
Докажем, что в этом случае $U(v)$ ---
непременно сектор. Во-первых, очевидно, что
$\ep \leq \alpha \leq \pi-\ep$,
и, значит, $U(v)$ --- не полоса и не полуплоскость.
Теперь предположим, что $U(v)$ --- двуугольник с углами $\alpha_1$ и $\alpha_2$.
Можно считать, что $\alpha_1\geq\alpha_2$.
Положим $\alpha(v)=\alpha_0$. Из определения
угловой функции $\alpha$  ясно, что $\alpha_0 \leq \alpha_1$.
 По теореме Гаусса--Бонне
$\alpha_1+\alpha_2=\omega^+(U(v))-\omega^-(U(v))$. Таким образом
$\alpha_1+\alpha_2 \leq \omega^+(U(v)) \leq \alpha_0 - \ep \leq
\alpha_1 - \ep < \alpha_1+\alpha_2$, что невозможно. Если же мы предположим, что
$U(v)$ --- дополнение к двуугольнику, то тогда,
как было отмечено выше, $\phi (-v)=-\phi(v)=(0,0)$, что,  как только что доказано,
невозможно
для двуугольника  $U(-v)$, что,,  невозможно.

Значит, мы получили разбиение полуплоскости на два сектора с
углами, удовлетворяющими нашим условиям. Теорема доказана.

\end{document}